\documentclass{amsart}
\usepackage{graphicx}
\usepackage{amsmath}
\usepackage{amssymb}
\theoremstyle{plain}

\newtheorem{thm}{Theorem}[section]
\newtheorem*{thm*}{Theorem}
\newtheorem{lem}[thm]{Lemma}
\newtheorem*{lem*}{Lemma}
\newtheorem{cor}[thm]{Corollary}
\newtheorem*{cor*}{Corollary}

\newtheorem*{cla*}{Claim}
\newtheorem{pro}[thm]{Proposition}
\newtheorem*{pro*}{Proposition}

\newtheorem*{fac*}{Fact}

\newtheorem*{que*}{Question}

\newtheorem*{con*}{Conjecture}
\newtheorem{rem}[thm]{Remark}
\newtheorem*{rem*}{Remark}

\newtheorem*{rems*}{Remarks}

\newtheorem*{defn*}{Definition}

\newcommand{\tr}{\mathrm{tr}}

\begin{document}

\title{Conjugate Generators of Knot and Link Groups}
\author{Jason Callahan}
\maketitle
\begin{abstract}
This note shows that if two elements of equal trace (e.g., conjugate elements) generate an arithmetic two-bridge knot or link group, then the elements are parabolic.  This includes the figure-eight knot and Whitehead link groups. Similarly, if two conjugate elements generate the trefoil knot group, then the elements are peripheral.
\end{abstract}

\section{Introduction}

By a knot or link group, we will mean the fundamental group of the knot or link complement in $S^3$.  It is well known that a two-bridge knot or link group is generated by two meridians of the knot or link (see \cite{BurdeZieschang03}).  The converse is also known; it is proved for link complements in $S^3$ in \cite{BoileauZimmermann89} (Corollary 3.3), for hyperbolic 3-manifolds of finite volume in \cite{Adams96} (Theorem 4.3), and for the most general class of 3-manifolds in \cite{BoileauWeidmann05} (Corollary 5):

\begin{thm} \label{converse}
If $M$ is a compact, orientable, irreducible 3-manifold with incompressible boundary and $\pi_1 M$ is generated by two peripheral elements, then $M$ is homeomorphic to the exterior of a two-bridge knot or link in $S^3$.
\end{thm}

Arising from work on Simon's Conjecture (see Section \ref{Simon} for statement and Problem 1.12 of \cite{Kirby97}), A.\ Reid and others proposed the following conjecture, which for convenience we will call:

\newtheorem*{rc}{Reid's Conjecture}
\begin{rc}
Let $K$ be a knot for which $\pi_1(S^3 \smallsetminus K)$ is generated by two conjugate elements. Then the elements are peripheral (and hence the knot is two-bridge by above).
\end{rc}

A knot in $S^3$ is a hyperbolic, satellite, or torus knot (Corollary 2.5 of \cite{Thurston82}).  By Proposition 17 of \cite{Gonzalez-AcunaRamirez02}, the $(p,q)$-torus knot group can be generated by two conjugate elements only when $p=2$, i.e., when the torus knot is two-bridge.  In Section \ref{torusknotsec}, we establish Reid's Conjecture for the $(2,3)$-torus knot group (i.e., the trefoil knot group), and we prove in Section \ref{Simon} that Reid's Conjecture implies Simon's Conjecture for two-bridge knots.  The latter conjecture has recently been proved in \cite{Boileauetal}.

When a knot or link complement in $S^3$ is hyperbolic, peripheral elements such as meridians map to parabolic elements under the discrete faithful representation of the knot or link group into $\mathrm{PSL}_2\mathbb{C}$.  Conjugate elements have equal trace, and we prove in Section \ref{figure8knot} a stronger version of Reid's Conjecture for the figure-eight knot (whose complement in $S^3$ is well known to be hyperbolic by \cite{Riley75}):

\begin{thm} \label{mainthm}
If two elements of equal trace generate the figure-eight knot group, then the elements are parabolic.
\end{thm}

The figure-eight knot group can, however, be generated by three conjugate loxodromic elements, so this result is, in some sense, sharp.

The proof of Theorem \ref{mainthm} relies heavily on the arithmeticity of the figure-eight knot complement in $S^3$.  Since the figure-eight knot is the only knot with arithmetic hyperbolic complement in $S^3$ (\cite{Reid91}), extending our result to all hyperbolic knot groups would require new techniques.  By Section 5 of \cite{GehringMaclachlanMartin98}, however, there are exactly four arithmetic Kleinian groups generated by two parabolic elements; each is the fundamental group of a hyperbolic two-bridge knot or link complement in $S^3$: the figure-eight knot, the Whitehead link, and the links $6^2_2$ and $6^2_3$.  In Section \ref{arithlinks}, we again exploit arithmeticity to extend our result to these:

\begin{thm} \label{otherthm}
If two elements of equal trace generate an arithmetic two-bridge knot or link group, then the elements are parabolic.
\end{thm}

\section{Preliminaries}

We collect several preliminary results that will be useful in what follows. The first is Lemma 7.1, together with the comments and definitions that precede it, in \cite{GehringMaclachlanMartinReid97}.

\begin{lem} \label{GMMR}
Let $\Gamma$ be a finite-covolume Kleinian group whose traces lie in $R$, the ring of integers in $\mathbb{Q}(\tr\Gamma)$.  If $\langle X, Y \rangle$ is a non-elementary subgroup of $\Gamma$, then $\mathcal{O} = R[1, X, Y, XY]$ is an order in the quaternion algebra
\[ A\Gamma = \left\{ \sum a_i \gamma_i : a_i \in \mathbb{Q}(\tr\Gamma), \gamma_i \in \Gamma \right\} \]
over $\mathbb{Q}(\tr\Gamma)$. Its discriminant $d(\mathcal{O})$ is the ideal $\langle 2 - \tr [X, Y] \rangle$ in $R$.
\end{lem}

The second is Theorem 6.3.4 in \cite{MaclachlanReid03}.

\begin{thm} \label{MR}
Let $\mathcal{O}_1$ and $\mathcal{O}_2$ be orders in a quaternion algebra over a number field.  If $\mathcal{O}_1 \subset \mathcal{O}_2$, then $d(\mathcal{O}_2) \, | \, d(\mathcal{O}_1)$, and $\mathcal{O}_1 = \mathcal{O}_2$ if and only if $d(\mathcal{O}_1) = d(\mathcal{O}_2)$.
\end{thm}

Our application is the following.

\begin{cor} \label{unitmult}
Let $\Gamma$ be a finite-covolume Kleinian group whose traces lie in $R$, the ring of integers in $\mathbb{Q}(\tr\Gamma)$.  If $\langle A, B \rangle = \Gamma = \langle X, Y \rangle$, then $2 - \tr [X,Y]$ is a unit multiple of $2 - \tr [A,B]$ in $R$.
\end{cor}

\begin{proof}
By Lemma \ref{GMMR}, $\mathcal{O}_1 = R[1, A, B, AB]$ and $\mathcal{O}_2 = R[1, X, Y, XY]$ are orders in $A\Gamma$.  Furthermore, $d(\mathcal{O}_1) = \langle 2- \tr[A, B] \rangle$ and $d(\mathcal{O}_2) = \langle 2 - \tr[X, Y] \rangle$ are ideals in $R$.  The Cayley-Hamilton Theorem yields the identity
\[ X + X^{-1} = \tr X \cdot 1 \]
which implies $A^{-1}, B^{-1} \in \mathcal{O}_1$ and $X^{-1}, Y^{-1} \in \mathcal{O}_2$.  Since $\mathcal{O}_1$ and $\mathcal{O}_2$ are ideals that are also rings with 1, we have
\[\Gamma = \langle A, B \rangle \subset \mathcal{O}_1 \text{ and } \Gamma = \langle X, Y \rangle \subset \mathcal{O}_2.\] Clearly, $R \subset \mathcal{O}_1,\mathcal{O}_2$, so $R\Gamma \subseteq \mathcal{O}_1,\mathcal{O}_2$, and $\mathcal{O}_1,\mathcal{O}_2 \subseteq R\Gamma$ by definition.  Therefore,
\[\mathcal{O}_1 = R\Gamma = \mathcal{O}_2.\]
Hence, $d(\mathcal{O}_1) = d(\mathcal{O}_2)$ by Theorem \ref{MR}, and the result follows.
\end{proof}

We conclude with two technical lemmas.

\begin{lem} \label{trace}
Let $\Gamma$ be a Kleinian group whose traces lie in $R$, the ring of integers in $\mathbb{Q}(\tr\Gamma)$, and $X, Y \in \Gamma$ with $\tr X = \tr Y$.  If
\[ x = \tr X = \tr Y, \; y = \tr XY - 2, \text{ and } z = 2 - \tr [X, Y], \]
then $y \, | \, z$ in $R$, and
\[ x^2 = \dfrac{z}{y} + y + 4. \]
\end{lem}

\begin{proof}
Standard trace relations (e.g., relation 3.15 in Section 3.4 of \cite{MaclachlanReid03}) yield
\begin{eqnarray*}
z &=& 4 + x^2\tr XY - \tr^2 XY - 2x^2 \\
  &=& (\tr XY - 2)x^2 - (\tr XY - 2)(\tr XY + 2) \\
  &=& y(x^2 - (y + 4))
\end{eqnarray*}
The result now follows.
\end{proof}

For the remainder of this note, let $\omega = e^{\pi i/3}$.

\begin{lem} \label{square}
Let $x = a + b\omega \in \mathbb{Z}[\omega]$, and
\[ x^2 = a^2 - b^2 + (2ab + b^2)\omega = n + m\omega \in \mathbb{Z}[\omega].\]
If $-4 \leq m \leq 4$, then the following are the only possibilities for $x^2$.
\begin{itemize}
  \item If $m = 0$, then $x^2 = a^2 \in \mathbb{Z}^2$ or $x^2 \leq 0$.
  \item If $m = \pm 1$, then $x^2 = -1 + \omega$ or $x^2 = -\omega$.
  \item If $m = \pm 2$, then $x^2 \not \in \mathbb{Z}[\omega]$, i.e., $m = \pm 2$ is not possible.
  \item If $m = \pm 3$, then $x^2 \in \{3\omega, -8 + 3\omega, 3 - 3\omega, -5 - 3\omega \}$.
  \item If $m = \pm 4$, then $x^2 = -4 + 4\omega$ or $x^2 = -4\omega$.
\end{itemize}
\end{lem}

\begin{proof}
For each case, we have the following.
\begin{itemize}
  \item If $m = 2ab + b^2 = 0$, then $b = 0$ or $a = -\dfrac{b}{2}$, so $x^2 = a^2 \in \mathbb{Z}^2$ or
  \[x^2 = -\dfrac{3b^2}{4} \leq 0.\]

  \item  If $m = 2ab + b^2 = \pm 1$, then $a = \dfrac{\pm 1 - b^2}{2b}$, so $b \, | \, 1$.  Thus, $b = \pm 1$, and
\[ (a,b) \in \{(0, \pm 1), \pm(1, -1) \}. \]
Therefore, $x^2 = a^2 - b^2 + (2ab + b^2)\omega = -1 + \omega$ or $x^2 = -\omega$.

  \item If $m = 2ab + b^2 = \pm 2$, then $a = \dfrac{\pm 2 - b^2}{2b}$, so $b \, | \, 2$ and $2 \, | \, b$.  Thus, $b = \pm 2$, and
\[ (a,b) \in \left\{\pm\left(\frac{1}{2}, -2\right), \pm\left(\frac{3}{2}, -2\right)\right\}. \]
Therefore, $x^2 = a^2 - b^2 + (2ab + b^2)\omega \not \in \mathbb{Z}[\omega]$.

  \item If $m = 2ab + b^2 = \pm 3$, then $a = \dfrac{\pm 3 - b^2}{2b}$, so $b \, | \, 3$.  Thus, $b \in \{\pm 1, \pm 3\}$, and
\[ (a,b) \in \{\pm(1, 1), \pm(1, -3), \pm(2, -1), \pm(2, -3) \}. \]
Therefore, $x^2 = a^2 - b^2 + (2ab + b^2)\omega \in \{3\omega, -8 + 3\omega, 3 - 3\omega, -5 - 3\omega \}$.

  \item If $m = 2ab + b^2 = \pm 4$, then $a = \dfrac{\pm 4 - b^2}{2b}$, so $b \, | \, 4$ and $2 \, | \, b$.  Thus, $b \in \{\pm 2, \pm 4\}$, and
\[ (a,b) \in \left\{\left(0, \pm 2\right), \pm\left(2, -2\right), \pm\left(\frac{3}{2}, -4\right), \pm\left(\frac{5}{2}, -4\right)\right\}.\]
Of these, the only values of $x^2 = a^2 - b^2 + (2ab + b^2)\omega \in \mathbb{Z}[\omega]$ are $x^2 = -4 + 4\omega$ and $x^2 = -4\omega$.
\end{itemize}
\end{proof}

\section{The Figure-Eight Knot} \label{figure8knot}

The discrete faithful representation of the figure-eight knot group
\[\pi_1 (S^3 \smallsetminus K) = \langle a,b \mid a^{-1}bab^{-1}a = ba^{-1}bab^{-1} \rangle\]
into $\mathrm{\mathrm{PSL}}_2\mathbb{C}$ is generated by
\[a \mapsto A = \left(\begin{array}{cc} 1&1 \\ 0&1 \\ \end{array}\right)
\text{ and }b \mapsto B = \left(\begin{array}{cc} 1&0 \\ \omega&1 \\
\end{array}\right)\]
(see \cite{Riley75}). Note our mild abuse of notation which will be continued tacitly throughout: we blur the distinction between elements of $\mathrm{PSL}_2\mathbb{C}$ and their lifts to $\mathrm{SL}_2\mathbb{C}$.  Then $\Gamma = \langle A, B \rangle$ is a finite-covolume Kleinian subgroup of $\mathrm{PSL}_2(\mathbb{Z}[\omega])$, so $\tr \Gamma \subset \mathbb{Z}[\omega]$, the ring of integers in the trace field $\mathbb{Q}(\omega)$. The invariant quaternion algebra is $M_2(\mathbb{Q}(\omega))$. Theorem \ref{mainthm} follows
immediately from:

\begin{thm}
If $\Gamma = \langle X, Y \rangle$ with $\tr X = \tr Y$, then $X$ and $Y$ are parabolic.
\end{thm}

\begin{proof}
By Corollary \ref{unitmult}, $2 - \tr[X, Y]$ is a unit multiple of $2 - \tr[A, B] = -\omega^2$ in $\mathbb{Z}[\omega]$.  The complete group of units in $\mathbb{Z}[\omega]$ is given by
\[(\mathbb{Z}[\omega])^* = \{1, \; \omega, \; \omega^2 = \omega - 1, \; \omega^3 = -1, \; \omega^4 = -\omega, \; \omega^5 = 1 - \omega\}.\]
Since $\Gamma$ is torsion-free, $\tr[X, Y] \not \in (-2,2)$, so
\[2 - \tr[X, Y] = \omega^n\text{ for some }n \in \{1, 2, 3, 4, 5\}.\]
Let \[ x = \tr X = \tr Y = a + b\omega \in \mathbb{Z}[\omega], \; y = \tr XY - 2, 
\text{ and } z = 2 - \tr [X,Y]. \]
Lemma \ref{trace} implies $y \, | \, z$ in $\mathbb{Z}[\omega]$, so $y$ is also a unit in $\mathbb{Z}[\omega]$.  Since $\tr XY \not \in (-2,2)$,
\[y = \omega^m \text{ for some }m \in \{0, 1, 2, 4, 5\}.\]
Varying $m$ and $n$ as above generates the following table of values for
\[x^2 = \dfrac{z}{y} + y + 4 = \omega^{n-m} + \omega^m + 4.\]

\[\begin{tabular}{l|l|l|l|l|l}
 &$ m=0 $&$ m=1 $&$ m=2 $&$ m=4 $&$ m=5 $\\
\hline
$n=1 $&$ 5 + \omega $&$ 5 + \omega $&$ 4 $&$ 3 - \omega $&$ 4 $\\
\hline
$n=2 $&$ 4 + \omega $&$ 4 + 2\omega $&$ 4 + \omega $&$ 4 - 2\omega $&$ 4 - \omega $\\
\hline
$n=3 $&$ 4 $&$ 3 + 2\omega $&$ 3 + 2\omega $&$ 5 - 2\omega $&$ 5 - 2\omega $\\
\hline
$n=4 $&$ 5 - \omega $&$ 3 + \omega $&$ 2 + 2\omega $&$ 5 - \omega $&$ 6 - 2\omega $\\
\hline
$n=5 $&$ 6 - \omega $&$ 4 $&$ 2 + \omega $&$ 4 $&$ 6 - \omega $\\
\end{tabular}\]

Of these, by Lemma \ref{square}, the only possible value for $x^2$
is 4, i.e., $X$ and $Y$ must be parabolic if they are to generate
$\Gamma$ and have equal trace.
\end{proof}

\begin{rem} \label{3conjlox}
The figure-eight knot group can be generated by three conjugate loxodromic elements.
\end{rem}

\begin{proof}
Let $\alpha = a^{-1}b^2$, $\beta = b\alpha b^{-1} = ba^{-1}b$, and 
$\gamma = b^{-1}\alpha b = b^{-1}a^{-1}b^3$. Then
\begin{eqnarray*}
\beta^{-1}\alpha\gamma^{-1}\alpha\beta^{-1}\alpha^2 &=& b^{-1}ab^{-1}a^{-1}b^2b^{-3}aba^{-1}b^2b^{-1}ab^{-1}a^{-1}b^2a^{-1}b^2 \\
&=& b^{-1}ab^{-1}a^{-1}b^{-1}a(ba^{-1}bab^{-1})a^{-1}b^2a^{-1}b^2 \\
&=& b^{-1}ab^{-1}a^{-1}b^{-1}a(a^{-1}bab^{-1}a)a^{-1}b^2a^{-1}b^2 \\
&=& b
\end{eqnarray*}
which implies $b \in \langle \alpha, \beta, \gamma \rangle$, so 
$b^2\alpha^{-1} = a \in \langle \alpha, \beta, \gamma \rangle$.  
Hence, $\langle \alpha, \beta, \gamma \rangle = \pi_1 (S^3 \smallsetminus K)$, and 
\[\alpha \mapsto \left(\begin{array}{cr} 1-2\omega&-1 \\ 2\omega&1 \\ \end{array}\right),
\beta \mapsto \left(\begin{array}{cc} 1-\omega&-1 \\ 1+\omega&1-\omega \\ \end{array}\right),
\text{ and }\gamma \mapsto \left(\begin{array}{cc} 1-3\omega&-1 \\ -3+5\omega&1+\omega \\
\end{array}\right).\]
Thus, the figure-eight knot group is generated three conjugate loxodromic elements.
\end{proof}

\section{The Arithmetic Two-Bridge Links} \label{arithlinks}

For the remainder of this note, let $\theta = \frac{1+i\sqrt{7}}{2}$.  The discrete faithful representation of
an arithmetic two-bridge link group into $\mathrm{PSL}_2\mathbb{C}$ is generated by the matrices
\[A = \left(\begin{array}{cc} 1&1 \\ 0&1 \\ \end{array}\right)
\text{ and }B = \left(\begin{array}{cc} 1&0 \\ \xi&1 \\ \end{array}\right)\]
where $\xi = 1+i$, $1+\omega$, and $\theta$ for the Whitehead link and the links $6^2_2$ and $6^2_3$ respectively (see Section 5 of \cite{GehringMaclachlanMartin98}). Then $\Gamma = \langle A, B \rangle$ is a finite-covolume Kleinian group whose traces lie in $\mathbb{Z}[i]$, $\mathbb{Z}[\omega]$, and $\mathbb{Z}[\theta]$ respectively.  These are the rings of integers in the respective trace fields $\mathbb{Q}(i)$, $\mathbb{Q}(\omega)$, and $\mathbb{Q}(\theta)$.  The respective invariant quaternion algebras are $M_2(\mathbb{Q}(i))$, $M_2(\mathbb{Q}(\omega))$, and $M_2(\mathbb{Q}(\theta))$. We now establish Theorem \ref{otherthm} for each arithmetic two-bridge link group.

\begin{thm}
Let $\Gamma$ be the discrete faithful representation in $\mathrm{PSL}_2(\mathbb{Z}[i])$ of the Whitehead link group with notation as above. If $\Gamma = \langle X, Y \rangle$ with $\tr X = \tr Y$, then $X$ and $Y$ are parabolic.
\end{thm}

\begin{proof}
By Corollary \ref{unitmult}, $2 - \tr[X, Y]$ is a unit multiple of $2 - \tr[A,
B] = -2i$ in $\mathbb{Z}[i]$.  Since $(\mathbb{Z}[i])^* = \{\pm 1, \pm i\}$
and $\tr[X, Y] \not \in (-2,2)$,
\[2 - \tr[X, Y] \in \{-2, \pm 2i\}.\]  Let
\[x = \tr X = \tr Y = a + bi \in \mathbb{Z}[i], \; y = \tr XY - 2, \text{ and } z = 2 - \tr[X, Y].\]
By Lemma \ref{trace}, $y \; | \; z$ in $\mathbb{Z}[i]$, so, since
$\tr XY \not \in (-2,2)$, we have
\[y \in \{1, \pm i, \pm (1+i), \pm (1-i), 2, \pm 2i\}.\] Varying
$y$ and $z$ as above generates the following table of possible
values for \[x^2 = \dfrac{z}{y} + y + 4.\]

\[\begin{tabular}{l|l|l|l}
 &$ z=-2 $&$ z=2i $&$ z=-2i $\\
\hline
$y=1 $&$ 3 $&$ 5+2i $&$ 5-2i $\\
\hline
$y=i $&$ 4+3i $&$ 6+i $&$ 2+i $\\
\hline
$y=-i $&$ 4-3i $&$ 2-i $&$ 6-i $\\
\hline
$y=1+i $&$ 4+2i $&$ 6+2i $&$ 4 $\\
\hline
$y=1-i $&$ 4-2i $&$ 4 $&$ 6-2i $\\
\hline
$y=-1+i $&$ 4+2i $&$ 4 $&$ 2+2i $\\
\hline
$y=-1-i $&$ 4-2i $&$ 2-2i $&$ 4 $\\
\hline
$y=2 $&$ 5 $&$ 6+i $&$ 6-i $\\
\hline
$y=2i $&$ 4+3i $&$ 5+2i $&$ 3+2i $\\
\hline
$y=-2i $&$ 4-3i $&$ 3-2i $&$ 5-2i $\\
\end{tabular}\]

We now check which of these values have the form $x^2 = a^2 - b^2 +
2abi \in \mathbb{Z}[i]$ based on 0, $\pm 1$, $\pm 2$, and $\pm 3$
being the only imaginary parts that arise in the table.

\begin{description}
\item[Case 1] The imaginary part of $x^2$ is 0; that is, \begin{eqnarray*}
  2ab = 0 &\Longrightarrow& a = 0\text{ or }b = 0 \\
  &\Longrightarrow& x^2 = -b^2\text{ or }x^2 = a^2 \in \mathbb{Z}^2
  \end{eqnarray*}
  The only value in the table with imaginary part 0 that can
  be expressed in either of these forms is $x^2 = 4$, i.e., $X$ and $Y$ are
  parabolic.
\item[Case 2] The imaginary part of $x^2$ is $\pm 1$; that is, $2ab = \pm
1$, which is impossible for $a, b \in \mathbb{Z}$.  Hence, $x^2$
cannot have imaginary part $\pm 1$.
\item[Case 3] The imaginary part of $x^2$ is $\pm 2$; that is, \begin{eqnarray*}
  2ab = \pm 2 &\Longrightarrow& a^2 = b^2 = 1 \\
  &\Longrightarrow& x^2 = \pm 2i
  \end{eqnarray*}
  But $\pm 2i$ does not appear in the table, so $x^2$ cannot have imaginary part $\pm 2$.
\item[Case 4] The imaginary part of $x^2$ is $\pm 3$; that is, $2ab = \pm
3$, which is impossible for $a, b \in \mathbb{Z}$.  Hence, $x^2$
cannot have imaginary part $\pm 3$.
\end{description}
This exhausts the table of possible values for $x^2$; therefore, $X$
and $Y$ must be parabolic if they are to generate $\Gamma$ and have
equal trace.
\end{proof}

\begin{thm}
Let $\Gamma$ be the discrete faithful representation in $\mathrm{PSL}_2(\mathbb{Z}[\omega])$ of the $6^2_2$ link group with notation as before. If $\Gamma = \langle X, Y \rangle$ with $\tr X = \tr Y$, then $X$ and $Y$ are parabolic.
\end{thm}

\begin{proof}
By Corollary \ref{unitmult}, $2 - \tr[X, Y]$ is a unit multiple of
$2 - \tr[A, B] = -3\omega$ in $\mathbb{Z}[\omega]$. Since
$\tr[X, Y] \not \in (-2,2)$, we have
\[2 - \tr[X, Y] \in \{\pm3\omega, \pm3(1-\omega), -3\} \subset 3(\mathbb{Z}[\omega])^*.\]
Let \[ x = \tr X = \tr Y = a + b\omega \in \mathbb{Z}[\omega], \; y
= \tr XY - 2, \text{ and } z = 2 - \tr [X,Y]. \] Lemma \ref{trace}
implies $y \, | \, z$, and hence $y \, | \, 3$, in
$\mathbb{Z}[\omega]$. Therefore, since $\tr XY \not \in (-2,2)$, we
have
\[y \in \{\pm(3-3\omega), \pm(2-\omega), \pm(1+\omega), \pm(1-\omega), \pm(1-2\omega), \pm\omega, \pm3\omega, 1, 3 \}.\]
Varying $y$ and $z$ as above generates the following table of values
for \[x^2 = \dfrac{z}{y} + y + 4.\]

\[\begin{tabular}{l|l|l|l|l|l}
 &$ z=3\omega $&$ z=-3\omega $&$ z=3-3\omega $&$ z=-3+3\omega $&$ z=-3 $\\
\hline
$y=3-3\omega $&$ 6-2\omega $&$ 8-4\omega $&$ 8-3\omega $&$ 6-3\omega $&$ 7-4\omega $\\
\hline
$y=-3+3\omega $&$ 2+2\omega $&$ 4\omega $&$ 3\omega $&$ 2+3\omega $&$ 1+4\omega $\\
\hline
$y=2-\omega $&$ 5+\omega $&$ 7-3\omega $&$ 8-2\omega $&$ 4 $&$ 5-2\omega $\\
\hline
$y=-2+\omega $&$ 3-\omega $&$ 1+3\omega $&$ 2\omega $&$ 4 $&$ 3+2\omega $\\
\hline
$y=1+\omega $&$ 6+2\omega $&$ 4 $&$ 6-\omega $&$ 4+3\omega $&$ 3+2\omega $\\
\hline
$y=-1-\omega $&$ 2-2\omega $&$ 4 $&$ 2+\omega $&$ 4-3\omega $&$ 5-2\omega $\\
\hline
$y=1-\omega $&$ 2+2\omega $&$ 8-4\omega $&$ 8-\omega $&$ 2-\omega $&$ 5-4\omega $\\
\hline
$y=-1+\omega $&$ 6-2\omega $&$ 4\omega $&$ \omega $&$ 6+\omega $&$ 3+4\omega $\\
\hline
$y=1-2\omega $&$ 3-\omega $&$ 7-3\omega $&$ 6-\omega $&$ 4-3\omega $&$ 6-4\omega $\\
\hline
$y=-1+2\omega $&$ 5+\omega $&$ 1+3\omega $&$ 2+\omega $&$ 4+3\omega $&$ 2+4\omega $\\
\hline
$y=\omega $&$ 7+\omega $&$ 1+\omega $&$ 4-2\omega $&$ 4+4\omega $&$ 1+4\omega $\\
\hline
$y=-\omega $&$ 1-\omega $&$ 7-\omega $&$ 4+2\omega $&$ 4-4\omega $&$ 7-4\omega $\\
\hline
$y=3\omega $&$ 5+3\omega $&$ 3+3\omega $&$ 4+2\omega $&$ 4+4\omega $&$ 3+4\omega $\\
\hline
$y=-3\omega $&$ 3-3\omega $&$ 5-3\omega $&$ 4-2\omega $&$ 4-4\omega $&$ 5-4\omega $\\
\hline
$y=1 $&$ 5+3\omega $&$ 5-3\omega $&$ 8-3\omega $&$ 2+3\omega $&$ 2 $\\
\hline
$y=3 $&$ 7+\omega $&$ 7-\omega $&$ 8-\omega $&$ 6+\omega $&$ 6 $\\
\end{tabular}\]

Of these, by Lemma \ref{square}, the only possible values for $x^2$
are 4, $3\omega$, and $3 - 3\omega$.  If
\[ x^2 = \tr^2 X = \tr^2 Y = 3\omega, \]
then $\tr X = \tr Y = \pm(1 + \omega)$, so the axes of $X$ and $Y$ project 
closed geodesics in $\mathbb{H}^3/\Gamma = S^3 \smallsetminus 6^2_2$ of length
\[ \text{Re}\left(2\cosh^{-1}\left(\pm\dfrac{1 + \omega}{2}\right)\right) \approx 1.087070145. \]
Similarly, if
\[ x^2 = \tr^2 X = \tr^2 Y = 3 - 3\omega, \]
then $\tr X = \tr Y = \pm(2 - \omega)$, so the axes of $X$ and $Y$ project 
closed geodesics in $\mathbb{H}^3/\Gamma = S^3 \smallsetminus 6^2_2$ of length
\[ \text{Re}\left(2\cosh^{-1}\left(\pm\dfrac{2 - \omega}{2}\right)\right) \approx 1.087070145. \]
But rigorous computation of the length spectrum in SnapPea 
(\cite{SnapPea}, \cite{HodgsonWeeks94}) shows that the shortest two closed
geodesics in $S^3 \smallsetminus 6^2_2$ have length 
\[0.86255462766206\text{ and }1.66288589105862.\]
Thus, the only possible value for $x^2$ is 4, i.e., $X$ and $Y$ must be parabolic 
if they are to generate $\Gamma$ and have equal trace.
\end{proof}

\begin{thm}
Let $\Gamma$ be the discrete faithful representation in $\mathrm{PSL}_2(\mathbb{Z}[\theta])$ of the $6^2_3$ link group with notation as before. If $\Gamma = \langle X, Y \rangle$ with $\tr X = \tr Y$, then $X$ and $Y$ are parabolic.
\end{thm}

\begin{proof}
By Corollary \ref{unitmult}, $2 - \tr[X, Y]$ is a unit multiple of
$2 - \tr[A, B] = -\theta^2 = 2 - \theta$ in
$\mathbb{Z}[\theta]$, i.e., $2 - \tr[X, Y] =\pm (2 - \theta)$.
Let \[ x = \tr X = \tr Y = a + b\theta \in \mathbb{Z}[\theta], \; y
= \tr XY - 2, \text{ and } z = 2 - \tr [X,Y]. \] Lemma \ref{trace}
implies $y \, | \, z$ in $\mathbb{Z}[\theta]$. Since
$\tr XY \not \in (-2,2)$, we have
\[y \in \{1, \pm (2 - \theta), \pm \theta\}.\]
Therefore,
\[x^2 = \dfrac{z}{y} + y + 4 \in \{4, 1+\theta, 3+\theta, 4 \pm 2\theta, 5-\theta, 7 - \theta\}.\]
Of these, arguing as before, the only possible value for
\[x^2 = a^2 - 2b^2 +(2ab + b^2)\theta \in \mathbb{Z}[\theta] \]
is 4, i.e., $X$ and $Y$ must be parabolic if they are to generate
$\Gamma$ and have equal trace.
\end{proof}

\section{Torus Knots} \label{torusknotsec}

Throughout this section, we assume $\gcd(p,q) = 1$ and $2 \leq p < q$. As is well known, the $(p,q)$-torus knot group admits the presentation
\[\pi_1 (S^3 \smallsetminus K_{p,q}) = \langle c,d \mid c^p=d^q\rangle,\]
which clearly surjects $\mathbb{Z}_p \ast \mathbb{Z}_q = \langle s,t \mid s^p = t^q = 1\rangle$ via $c \mapsto s,$ $d \mapsto t$ (see \cite{BurdeZieschang03}).

We begin our investigation of conjugate generators for torus knot groups by paraphrasing Proposition 17 of \cite{Gonzalez-AcunaRamirez02}:

\begin{pro} \label{GAR}
If $\mathbb{Z}_p \ast \mathbb{Z}_q = \langle s,t \mid s^p = t^q = 1 \rangle$ can be generated by two conjugate elements, then $p = 2$.
\end{pro}

Thus, via the surjection above, the $(p,q)$-torus knot group can be generated by two conjugate elements only when $p=2$, i.e., when the torus knot is two-bridge with normal form $(q/1)$. The $(2,q)$-torus knot group $\langle c,d \mid c^2=d^q \rangle$ has a parabolic representation into $\mathrm{PSL}_2\mathbb{C}$ generated by
\[c \mapsto C = \left(\begin{array}{cc} 0&(2\cos(\pi/q))^{-1} \\ -2\cos(\pi/q)&0 \\ \end{array}\right)
\text{ and }d \mapsto D = \left(\begin{array}{cc} 1-4\cos^2(\pi/q)&1 \\ -4\cos^2(\pi/q)&1 \\ \end{array}\right)\]
(Theorem 6 of \cite{Riley72}). Then $\Gamma_q = \langle C,D \rangle$ is a finite-coarea subgroup of $\mathrm{PSL}_2\mathbb{R}$ and has a presentation $\langle C,D \mid C^2 = D^q = 1\rangle = \mathbb{Z}_2 \ast \mathbb{Z}_q$.  Furthermore, the traces of $\Gamma_q$ are algebraic integers in $\mathbb{Q}(\tr \Gamma) = \mathbb{Q}(\cos(\pi/q))$ (cf. Section 8.3 of \cite{MaclachlanReid03}).

To analyze the case $q=3$ (i.e., the trefoil knot), we recall the classification of generating pairs for the groups $\mathbb{Z}_p \ast \mathbb{Z}_q$ up to Nielsen equivalence (Corollary 4.14 of \cite{Zieschang77} and Theorem 13 of \cite{Gonzalez-AcunaRamirez02}):

\begin{thm} \label{Z}
Every generating pair for the group $\mathbb{Z}_p \ast \mathbb{Z}_q = \langle s,t \mid s^p = t^q = 1\rangle$ is Nielsen equivalent to exactly one generating pair of the form $(s^m,t^n)$ where 
\[\gcd(m,p)=\gcd(n,q)=1, \; 0<2m\leq p, \; and \; 0<2n\leq q.\]
\end{thm}

\begin{cor} \label{modular}
If two conjugate elements generate the $(2,3)$-torus knot group (i.e., the trefoil knot group), then the elements are peripheral.
\end{cor}

\begin{proof}
Suppose $\langle \alpha, \beta \rangle = \pi_1 (S^3 \smallsetminus K_{2,3})$ with $\alpha$ conjugate to $\beta$. Let $\alpha \mapsto X$ and $\beta \mapsto Y$ via the representation above. Then $\langle X,Y \rangle = \Gamma_3 = \mathbb{Z}_2 \ast \mathbb{Z}_3$, and $X$ is conjugate to $Y$, so $\tr X = \tr Y$. Thus, by Theorem \ref{Z} with $p=2$ and $q=3$, $(X,Y)$ is Nielsen equivalent to $(C,D)$. Since commutators of Nielsen equivalent pairs have equal trace, we have
\[2 - \tr[X,Y] = 2 - \tr[C,D] = -1.\]
Lemma \ref{trace} then implies $\tr XY - 2 = \pm 1$, and so $\tr^2 X = \tr^2 Y = 4$. Hence, $\alpha$ and $\beta$ are peripheral.
\end{proof}

As defined in this section, $\Gamma_q$ is a discrete faithful representation into $\mathrm{PSL}_2\mathbb{R}$ of the $(2,q,\infty)$-triangle group. Thus, the proof of Corollary \ref{modular} shows that if two elements of equal trace (e.g., conjugate elements) generate the $(2,3,\infty)$-triangle group (i.e., the modular group), then the elements are parabolic.  Following Section 13.3 of \cite{MaclachlanReid03} (cf. \cite{Takeuchi77}), the $(2,q,\infty)$-triangle group is arithmetic only when $q = 3$, 4, 6, or $\infty$. Similar methods can then be used to show that if two elements of equal trace generate an arithmetic $(2,q,\infty)$-triangle group, then the elements are parabolic.  The $(2,q,\infty)$-triangle group can, however, be generated by two conjugate hyperbolic elements when $q>3$ is odd: with $\Gamma_q = \langle C,D \rangle$ as above, let $X = CD$ and $Y = C^{-1}XC = DC$; then $(YX)^{\frac{q+1}{2}} = D$, so $\langle X,Y \rangle = \Gamma_q$, and $\tr X = \tr Y = -4\cos(\pi/q) < -2$ since $q>3$, so $X$ and $Y$ are hyperbolic.

\section{Simon's Conjecture} \label{Simon}

The following is attributed to J.\ Simon (cf. Problem 1.12 of \cite{Kirby97}).

\newtheorem*{simon}{Simon's Conjecture}
\begin{simon}
A knot group can surject only finitely many other knot groups.
\end{simon}

To show Reid's Conjecture implies Simon's Conjecture for two-bridge knots, we first note Theorem 5.2 of \cite{ReidWang99} (recall that a knot complement is called {\em small} if it does not contain a closed embedded essential surface):

\begin{thm} \label{reidwang}
If $M$ is a small hyperbolic knot complement in $S^3$, then there exist only finitely many hyperbolic 3-manifolds $N$ for which there is a peripheral preserving epimorphism $\pi_1 M \twoheadrightarrow \pi_ 1 N$.
\end{thm}

\begin{thm}
Reid's Conjecture implies Simon's Conjecture for two-bridge knots.
\end{thm}

\begin{proof}
Let $K$ be a two-bridge knot and 
$\varphi_i : \pi_1 (S^3 \smallsetminus K) \twoheadrightarrow \pi_1 (S^3 \smallsetminus K_i)$
a collection of epimorphisms $\varphi_i$ and knots $K_i$.

If $K_i$ is a $(p,q)$-torus knot, then its Alexander polynomial has degree $(p-1)(q-1)$ and divides the Alexander polynomial of $K$ since $\varphi_i$ is an epimorphism (see \cite{BurdeZieschang03}).  This can occur for only finitely many integer pairs $(p,q)$; hence, only finitely many $K_i$ are torus knots.

Now assume $K_i$ is not a torus knot.  Since $K$ is two-bridge, $\pi_1 (S^3 \smallsetminus K)$ is generated by two conjugate meridians, $a$ and $b$.  Then $\varphi_i (a)$ and $\varphi_i (b)$ are conjugate elements that generate $\pi_1 (S^3 \smallsetminus K_i)$, so Reid's Conjecture implies that they are peripheral.  Hence, by Theorem \ref{converse}, $K_i$ is two-bridge and therefore hyperbolic since it is not a torus knot.

Let $\lambda \in \pi_1 (S^3 \smallsetminus K)$ such that $\langle a, \lambda \rangle$ is a peripheral subgroup of $\pi_1 (S^3 \smallsetminus K)$.  Then $\lambda$ commutes with $a$ in $\pi_1 (S^3 \smallsetminus K)$, so $\varphi_i(\lambda)$ commutes with $\varphi_i(a)$ in $\pi_1 (S^3 \smallsetminus K_i)$ and hence is a peripheral element of $\pi_1 (S^3 \smallsetminus K_i)$. Therefore, $\varphi_i$ is peripheral preserving, so the result follows from Theorem \ref{reidwang} since hyperbolic two-bridge knot complements in $S^3$ are small (\cite{HatcherThurston85}).
\end{proof}

\subsection*{Acknowledgements}
The author thanks his advisor, Alan Reid, whose patience and guidance made this work possible, and Eric Chesebro for assistance with Remark \ref{3conjlox}.

\bigskip

\address{
\noindent The University of Texas at Austin\\
Department of Mathematics\\
1 University Station C1200\\
Austin, TX 78712, USA}

\medskip

\noindent\email{callahan@math.utexas.edu}

\end{document}